\documentclass[12pt]{iopart}
\usepackage{iopams}

\usepackage{graphicx}

\def\R{{\mathbb R}}

\def\Z{{\mathbb Z}}

\def\ds{\displaystyle}

\def\all{\forall}

\def\sl2{SL(2,\R)}

\def\csl2{\widetilde{\sl2}}

\def\x{{\bf x}}
\def\y{{\bf y}}
\def\z{{\bf z}}
\newcommand{\bea}{\begin{eqnarray}}
\newcommand{\eea}{\end{eqnarray}}
\newcommand{\beq}{\begin{equation}}
\newcommand{\eeq}{\end{equation}}
\newcommand{\barr}{\begin{array}}
\newcommand{\earr}{\end{array}}
\newcommand{\ovl}{\overline}
\newtheorem{theo}{Theorem}[section]
\newtheorem{prop}{Proposition}[section]
\newtheorem{defn}{Definition}[section]

\newtheorem{cor}{Corollary} [section]

\newtheorem{rmk}{Remark}[section]
\def\bth{\begin{theo}}
\def\eth{\end{theo}}
\def\bpr{\begin{prop}}
\def\epr{\end{prop}}
\def\bdf{\begin{defn}}
\def\edf{\end{defn}}
\def\brmk{\begin{rmk}}
\def\ermk{\end{rmk}}

\begin{document}

\title{Monotone gradient dynamics  and location of stationary $(p,q)$--configurations}
\author{Emilia Petrisor}


\begin{abstract}  Exploiting the monotone property of the gradient dynamics of  the Frenkel--Kontorova model,  we locate in the space of $(p,q)$--configurations the ordered and unordered stationary states, as well as forbidden regions for such states. Moreover we show that some generalized Frenkel--Kontorova models
(associated to  multiharmonic standard maps) can have
ordered $(p,q)$--configurations  that are neither action minimizing nor minimaximizing, and give their location with respect to the set of $(p,q)$--minimizers and minimaximizers.
\end{abstract}
\ams{37E40}
\section{Introduction} Area preserving twist maps have been extensively studied  as being typical examples of
 dynamical systems which exhibit a full range of  behaviour, from regular to chaotic motions.
Their  study  is mainly concerned with the existence/breakup of rotational invariant circles (RICs). These invariant sets are barrier to transport and their breakup leads to loss of stability.   A criterion that gives  conditions ensuring that no RIC of given rotation number can exist or no RIC can pass through a point or region of the phase space is called converse KAM criterion \cite{kayperciv}.
  A survey of the known converse KAM criteria to date
can be found in \cite{calejallave}.\par

One of the first such  criteria was formulated by
Boyland--Hall  \cite{boyland1}, \cite{boyland2}. It states that if a twist map has an unordered
periodic orbit, $\mathcal{O}$, then there are no invariant circles whose
 rotation number  belongs to the  rotation band of $\mathcal{O}$.\par
 The search for unordered periodic orbits using  classical methods can be difficult due to the instability of the system in the region of the phase space where they lie.
  It is much more  simple  to detect them looking for their representatives in the space of $(p,q)$--configurations,
i.e. as    zeros of  the periodic action gradient.\par

The gradient of the  action functional of an area preserving twist map  of the infinite annulus (cylinder)  is in fact the gradient of the energy  of the Frenkel--Kontorova (FK) model, studied in solid--state physics \cite{FloriaA}, \cite{FloriaB}.
The standard FK model  is associated to an infinite one-dimensional  chain of atoms connected by harmonic springs and subjected to an external  potential \cite{aubry}, \cite{FloriaB}. If the potential is multiharmonic, then the corresponding FK model is a also called generalized FK model.\par

Using an interplay between the variational approach to the dynamics of area preserving twist maps, and the gradient dynamics of the periodic action  we locate the regions in the space of $(p,q)$--configurations where ordered or unordered stationary states of the action gradient can lie, respectively the forbidden regions for such states.
More precisely, we show that  if
$\x, \x'$,  are  two consecutive $(p,q)$--minimizers, and  $\y$ is a $(p,q)$--minimaximizer,  such that $\x\prec \y\prec {\x}'$, then no other stationary $(p,q)$--configuration can belong to the intervals $[[\x,\y]]$, $[[\y, \x']]$ (Proposition \ref{nostatincon}), but in the complement of
the set $[\x, \y]\cup [\y, \x']$, with respect to  the interval $[[x,x']]$ there can exist ordered  $(p,q)$--stationary configurations.  We give an example of three  harmonic standard map, exhibiting an ordered  $(1,2)$--orbit which is neither minimizing nor minimaximizing, and whose corresponding configuration lies in such a complement.  As far as we know no such periodic orbit was revealed  in the
  theoretical and numerical investigations of the twist maps, reported in the literature.

\par Location of the $(p,q)$--stationary states
has both a theoretical and practical importance. In an attempt to implement the Hall--Boyland criterion, first of all we need to know where in the space of $(p,q)$--configurations it is advisable to look for stationary states, and then to decide whether an identified  such a configuration is ordered or not.\par\medskip

 A   general presentation  of the gradient dynamics of the Frenkel--Kontorova model can be found in \cite{FloriaA}, \cite{baesens}. Angenent \cite{angenent} was the first who exploited the monotonicity of this gradient semiflow in the study  of periodic orbits of twist maps. Gol\'{e} \cite{gole2001},  Mramor\&Rinky \cite{mramor} and Slijep\v{c}evi\'{c} \cite{slijepcevic} have also used the monotonicity of the gradient semiflow of the FK--model to prove the existence of ghost circles, respectively to prove
 the existence of Mather's shadowing orbits of twist maps.\par\medskip
We consider  an exact $C^r$--area preserving positive twist diffeomorphism ($r\geq 1$), $f$,  of the infinite annulus  $\mathbb{A}=\mathbb{T}^1\times \R$, represented by a lift, $F$  (for basic properties of such maps the reader is referred to \cite{gole2001}). Denote by $\pi:\R^2\to\mathbb{A}$ the covering projection.\par
Dynamics of $F$ has a variational formulation. Its
 orbits are in a one--to-one correspondence with the critical
points of an action. \par
 Under the hypotheses on $f$ it follows that the lift $F$ admits a  $C^2$--generating function  $h:\R^2\to\R$ (unique up to additive constants), such that
$h_{12}(x,x')=\partial h/\partial x\partial x'<0$ on $\R^2$, and
 $F(x,y)=(x', y')$ iff $y=-h_1(x,x')$, $y'=h_2(x,x')$  ($h_1=\partial h/\partial x, h_2=\partial h/\partial x')$.\par
  For $p,q$  coprime integers, $q>0$, the $F$--orbit of a point $(x_0,y_0)\in\R^2$ is called a $(p,q)$--type orbit,  and its projection onto annulus  a $(p,q)$-periodic orbit if $F^q(x_0,y_0)=(x_0+p,y_0)$.\par
  We denote by $\mathcal{X}_{pq}$, the space of $(p,q)$--type configurations, i.e. sequences  of real numbers, ${\bf x}=(x_n)_{n\in\Z}$,
 such that $x_{n+q}=x_n+p$, for all
$n\in\Z$.   Being  an affine
subspace of $\R^{q+1}=\{(x_0, x_1, \ldots, x_q)\}$, of equation
$x_q=x_0+p$, $\mathcal{X}_{pq}$ can be identified with $\R^q$.\par

The generating function $h$ defines the action
$W_{pq}$, on the space $\mathcal{X}_{pq}$ of $(p,q)$--configurations:
  $$
W_{pq}({\bf x})=\sum_{k=0}^{q-1}h(x_k, x_{k+1})$$
By Aubry--Mather theory, \cite{aubry},\cite{mather}, for each  pair of relative prime integers
$(p,q)$, $q>0$, $F$  has at least two orbits  of
type $(p,q)$. One corresponds to a non--degenerate  minimizing
$(p,q)$--configuration of the action.   Associated
to such an orbit there is a second one, corresponding to a
mini--maximizing configuration.  The corresponding $F$--orbits and $f$--orbits are called in the sequel $(p,q)$--minimizing, respectively $(p,q)$--minimaximizing orbits.  These orbits are well ordered.\par  More precisely,
an invariant set $M$ of a positive twist map, $f$, is well
ordered if  for every $(x,y), (x', y')'\in\pi^{-1}(M)$, we have $x<x'$ iff
$F_1(x,y)<F_1(x',y')$, where $F_1$ is the first component of the
lift $F$.

A well--ordered $(p,q)$--orbit is also called Birkhoff orbit, while
a badly--ordered orbit is called non--Birkhoff or unordered orbit.\par

An appropriate framework to study
 the order properties of periodic
orbits is defined as follows. \par The space $\mathcal{X}_{pq}$ of $(p,q)$--configurations is
partially ordered with respect to an order relation inherited from $\R^\Z$.
   $\x=(x_k)$, $\y=(y_k)\in\R^\Z$ are related, and we write:
\begin{equation}\label{seqorder}\barr{ll} \x \leq \y&\Leftrightarrow\quad
x_k\leq y_k, \all\,\, k\in\Z\earr\end{equation}

One also defines: \begin{equation}\label{relless}\barr{ll}\x
<\y&\Leftrightarrow\quad \x \leq \y, \mbox{but}\,\, \x\neq
\y\\
\x\prec \y&\Leftrightarrow\quad  x_k<y_k, \all\,\, k\in\Z\earr \end{equation}
If $\x<\y$ one says that $\x$ and $\y$ are {\it weakly ordered}, while if $\x\prec \y$ one says that they are {\it strictly ordered}.
$(\mathcal{X}_{pq},\leq)$ is a lattice.\par
To each configuration $\x\in  \mathcal{X}_{pq}$ one associates the positive order cone
$V_+(\x)=\{{\bf y}\in\mathcal{X}_{pq}\,|\, \x\leq \y\}$, respectively the negative
order cone $V_-(\x)=\{{\bf z}\in\mathcal{X}_{pq}\,|\, {\bf z}\leq \x\}$. Any  two  $(p,q)$--configurations, $\x, \y\in \mathcal{X}_{pq}$, comparable with respect to the relation  $\leq$, or  $\prec$, define  respectively the intervals:
 $$[\x, \y]=\{{\bf z}\in\mathcal{X}_{pq}\,|\, \x\leq
{\bf z}\leq \y\}, \quad  [[\x,\y]]=\{{\bf z}\in\mathcal{X}_{pq}\:|\: x\prec {\bf z}\prec \y\}$$

Let $\tau_{ij}:\R^\Z\to\R^\Z$ be the translation
map defined by:
$$(\tau_{ij}\x)_k=x_{k+i}+j, \quad \all\,\, \x=(x_k)\in\R^\Z, i,j,k\in\Z$$
$\mathcal{X}_{pq}$ is invariant to any integer translation, $\tau_{ij}$.\par
 A $(p,q)$--configuration, $\x$,  such that:
\begin{equation}\label{COorder}\all\,\, i,j\in\Z,  \:\mbox{either}\:  \x\leq \tau_{ij}\x\quad \mbox{or}\quad
\tau_{ij}\x\leq\x\end{equation}
is called cyclically ordered, as well as the corresponding $(p,q)$--orbit of $F$:
\begin{equation}\label{criticseqpq} (x_k, y_k)=(x_k, -h_1(x_k, x_{k+1})),  \quad \all\,\, k\in\Z\end{equation}

Let $\x$ be a   $(p,q)$--configuration. The piecewise affine function  that interpolates linearly
the points  $(k, x_k), k\in\Z$ is  called Aubry function.
  \par
Let $M_{pq}$ be the subset of $(\mathcal{X}_{pq},\leq)$ consisting
in all  $W_{pq}$--globally  minimizing configurations. Any  $\x\in M_{pq}$ is  cyclically ordered. By Aubry--Mather theory \cite{mather} $M_{pq}$ is a completely ordered subset $\mathcal{X}_{pq}$. In the sequel  the elements in $M_{pq}$ are referred to as $(p,q)$--minimizers. \par  $\x, \x'\in M_{pq}$ are called consecutive if $\x\prec \x'$ and there is no ${\bf x}''\in M_{pq}$ such that  $\x\prec {\bf x}''\prec \x'$.\par

 \par\medskip

Besides the   variational  interpretation, the $(p,q)$-- orbits of the twist map $F$ can be also associated to the stationary states (equilibrium points) of the gradient of the action, $\nabla W_{pq}$. If the second derivative of the generating function
$h$ of  the  twist map $F$ is bounded, then the system of differential equations
$\dot{\x}=-\nabla W_{pq}$    is a cooperative dynamical system  on $(\mathcal{X}_{pq},\leq)$ (in fact $\mathcal{X}_{pq}$ is invariant to the semiflow of the minus gradient of the energy of the Frenkel--Kontorova model; for details see  \cite{gole2001}, \cite{baesens}).
Its   $C^1$--semiflow, $\xi_t$, $t\geq 0$,  commutes with the group of translations $\{\tau_{ij}\:|\: i,j, \in\Z\}$. Any stationary state, $\x\in\mathcal{X}_{pq}$, is a critical $(p,q)$--configuration  of the action and  the linearized gradient at $\x$, $D^2W_{pq}(\x)$, is the Hessian of $W_{pq}$ at $\x$.\par

An important property of the gradient of the action relevant for our approach is that the semiflow of the minus gradient   is strictly monotone (\cite{gole2001}, \cite{baesens}), i.e.
\begin{equation}\label{strongMonon}
\x<\y\Rightarrow \xi_t\x\prec \xi_t\y, \all\,\, t>0\end{equation}

With these results in mind, in the next section we detect the  regions in the space of $(p,q)$--configurations where no stationary state  can
lie, respectively where  cyclically ordered or unordered $(p,q)$--stationary states is possible to be located.\par

\section{Location of stationary $(p,q)$--configurations}
We denote by $\mathcal{S}_{pq}=\{\z\in\mathcal{X}_{pq}\:|\: \nabla W_{pq}(\z)=0\}$ the set of stationary states in the space of $(p,q)$--configurations.
\bpr\label{nostateBcone}  Let  $\x$ be a minimizer or a minimaximizer of  $W_{pq}$. Then no stationary state $\z\in \mathcal{X}_{pq}$, $z\neq \x$, can belong to the boundaries of the cones
$V_{\pm}(\x)$.\epr
Proof:  Suppose $\z\in \mathcal{S}_{pq}$,  $\z\neq \x$, and $\z\in \partial V_+(\x)$ (boundary of $V_+(\x)$). This means that  $\x <\z$ and from the strict monotonicity of the flow  it follows that $\xi_t\x\prec \xi_t\z$. But the last relation is impossible because $\xi_t\x=\x$, $\xi_t\z=\z$, $\all t\geq 0$. Similarly one shows that $\z\notin \partial V_-(\x)$. \par\medskip

\bpr\label{regordered}  If  $\z \in\mathcal{S}_{pq}\setminus M_{pq}$ then there exist   two consecutive $(p,q)$--minimizers, $\x\prec \x'$, such that  $x_0<z_0<x'_0$. If $\z \in \mbox{Int}(V_+(\x)\cap V_-(\x'))=[[\x, \x']]$, then $\z$ is cyclically ordered, while if  $\z\in \mbox{Int}(V_+(\x))\setminus V_-(\x')$ or $\z\in \mbox{Int}(V_-(\x'))\setminus V_+(\x)$, then $\z$ is unordered.\epr

Proof:  Since to each $\x\in \mathcal{S}_{pq}$ one associates the $F$--orbit of the point $(x_0, -h_1(x_0, x_1))\in\R^2$, it is obvious the existence of the two consecutive minimizers as stated.\par
If  $\z\in [[\x, \x']]$, then for any $i, j\in\Z$, $\tau_{ij}\x\prec \tau_{ij}\z\prec \tau_{ij}\x'$. Since $\x, \x'$ are consecutive elements in $M_{pq}$, it follows that no translation , $\tau_{k\ell}$, can exist, such that  $\tau_{k\ell}\x\in [[\x, \x']]$ or $\tau_{k\ell}\x'\in [[\x, \x']]$. Thus
 except for $i=0, j=0$ the intervals $[[\x, \x']]$, $[[\tau_{ij}\x, \tau_{ij}\x']]$ have no point in common. Hence  $\x\leq \tau_{ij}\x$
implies $\z\leq \tau_{ij}\z$, and analogously $\tau_{ij}\x\leq \x$ implies $\tau_{ij}\z\leq \z$. Since $\x$ is cyclically ordered it follows that $z$ is also cyclically ordered. One also says that $\z$ is cyclically ordered with respect to $M_{pq}$.\par
If $\z\notin [\x, \x']$, but $x_0<z_0<x'_0$, then either   $\z-\x$ or $\z-\x'$ does not have coordinates
of the same sign. Since the Aubry functions associated to  $\x$ and $\x'$  are monotone, it follows that the Aubry function of $\z$ is not monotone, and as a consequence $\z$ cannot be cyclically ordered.

\begin{cor} A   stationary state $y\in\mathcal{X}_{pq}$ which is ordered with respect to $M_{pq}$ is comparable with each $x\in M_{pq}$, while an unordered one is incomparable with at least a  $(p,q)$--minimizer.\end{cor}
\begin{rmk} From Aubry--Mather theory we know that the $(p,q)$-minimaximizing configurations
are
 cyclically ordered with respect to elements in $M_{pq}$ \cite{mather}.  More precisely,  if the action $W_{pq}$  is a Morse function (all its critical points are nondegenerate), then for any  consecutive $(p,q)$--minimizers, $\x\prec {\x}'$, there exists a minimaximizing sequence $\y\in\mathcal{X}_{pq}$ such that
$\x\prec \y\prec {\x}'$ and from the above Proposition it follows that it is cyclically ordered with respect to $M_{pq}$. \end{rmk}

 We are wondering whether besides the $(p,q)$--minimaximizers
there can exist other cyclically ordered $p,q)$--stationary states with respect to $M_{pq}$.
In order to give an answer  let us first analyze where a cyclically ordered stationary state, different from a $(p,q)$--minimaximizer can lie within an interval $[[\x, \x']]$, of ends
$\x, \x'$ that are  consecutive $(p,q)$--minimizers.
\bpr\label{nostatincon}  If   $\x,\x'$  are  two consecutive $(p,q)$--minimizers and $\y\in \mathcal{X}_{pq}$ a minimaximizer,  such that $\x\prec \y\prec {\x}'$, then no stationary state of the gradient semiflow, $\xi_t$, can belong to the intersection of order cones  $V_+(\x)\cap V_-(\y)$, $V_+(\y)\cap V_-({\x}')$.\epr
Proof:
The derivative $-D^2 W_{pq}$  at $\x'$ has $q$ negative eigenvalues, while at $\y$, one positive eigenvalue. Thus the unstable manifold of $\y$ is one--dimensional. Moreover, it is
 a strictly ordered  heteroclinic connection, $\gamma$, between the two hyperbolic equilibria (an arc of ghost circle) \cite{gole2001},  \cite{mramor}  \cite{fuscoliva}) (a detailed proof of the existence of this heteroclinic connection is given in \cite{mramor}, Lemma 8.7). For every ${\bf v}\in \gamma$, $\xi_t {\bf v}\to \x'$, as $t\to\infty$.\par
Suppose that that there is a $(p,q)$--stationary state $\z\in V_+(\y)\cap V_-({\x}')$. Since $\z$ cannot belong to the boundaries of the two cones, it is  interior to $V_+(y)\cap V_-(x')$, and we can choose a point ${\bf v}\in \gamma$ such that ${\bf v}<{\bf z}$.
It follows that  $\xi_t {\bf v}\prec \xi_t {\bf z}$, $\all\, \, t>0$, which is impossible because $\xi_t {\bf v}\to \x'$, as $t\to\infty$, and $\xi_t \z=\z$, $\all\,\, t\geq 0$. Similarly one shows that $\z$ cannot belong to the intersection $V_+(\x)\cap V_-(\y)$.\par\medskip

\begin{figure}[h]
\begin{center}
\includegraphics{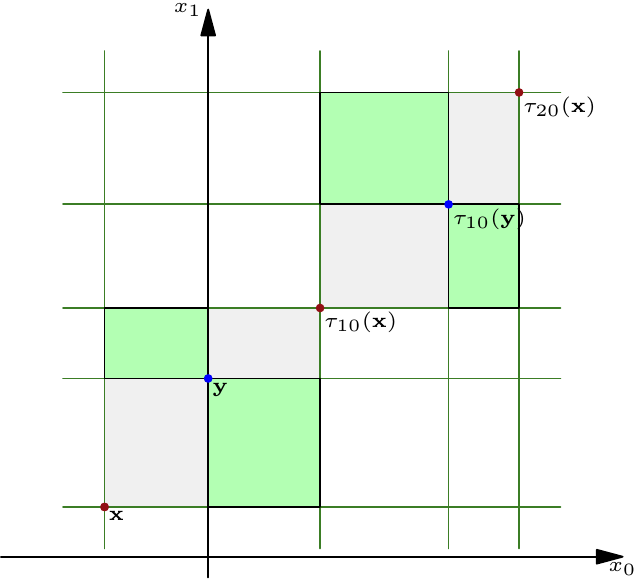}
\caption{\label{stdLaticeaXpqN}  The space of $(1,2)$--configurations. Details in text.}
\end{center}
\end{figure}

In Fig.\ref{stdLaticeaXpqN}  is illustrated a part of the space of $(1,2)$ configurations. $\x$ and its integer translates $\tau_{i0}$, $i=1,2$, are $(1,2)$--minimizers, $\y$ and $\tau_{10}(\y)$ are minimaximizers.  The gray semilines starting from these points are boundaries for the associated order cones. The semilines as well as the gray regions cannot contain other stationary $(1,2)$--configurations. In green regions, theoretically on can find ordered $(1,2)$--stationary states, while in the white regions, unordered ones.

\begin{figure}[h]
\begin{center}
\includegraphics{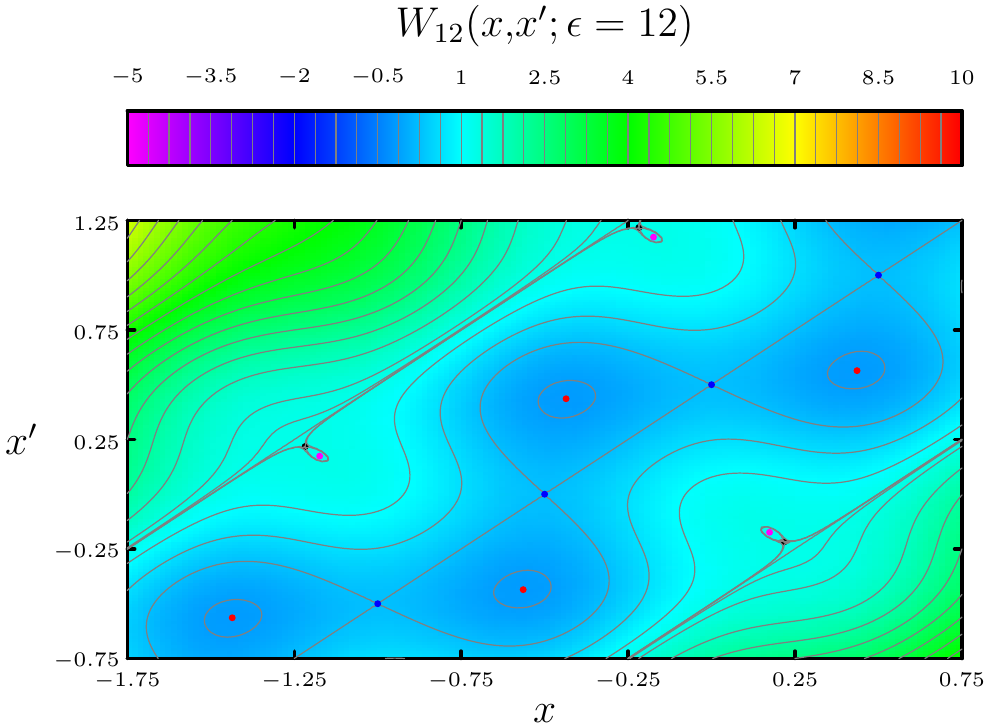}
\caption{\label{stdcontWpqk12}  The contour plot of the action $W_{12}$ defined by the generating function of the standard map, corresponding to $\epsilon=12$, and ordered and unordered $(1,2)$--stationary configurations. }
\end{center}
\end{figure}

  In Fig.\ref{stdcontWpqk12} is given the contour plot of the action $W_{12}$, associated to the generating function $h(x,x')=\ds\frac{1}{2}(x-x')^2+\ds\frac{\epsilon}{(2\pi)^2}\cos(2\pi x)$ of the standard map. The red points are  $(1,2)$--minimizers, those colored in blue are $(1,2)$--minimaximizers, while lateral points are unordered $(1,2)$--stationary configurations. One can observe that all these points have relative position as we have deduced above.

We note that  in $\mathcal{X}_{pq}/\tau_{0,1}:=\mathcal{X}_{pq}/G$ ($G$ is a group isomorphic with $\Z$, generated by the translation $\tau_{01}$) there exist at least $q$ minimizers. If a positive  twist map exhibits only a $(p,q)$--minimizing orbit (this is a case of the classical standard map)
then in $\mathcal{X}_{pq}/\tau_{0,1}$ there exist exactly $q$ minimizers. Actually if
$\x=(x_0, x_1, \ldots, x_{q-1})$ is a $(p,q)$--minimizer then there exist $q$ pairs  $(i,j)\in\Z\times \Z$ such that  $(\tau_{ij}\x)_0\in[0,1)$.
 Fig.\ref{stdLaticeaXpqN} is representative for such a case.\par

 Computer experiments reveal  that a multiharmonic standard map, i.e. a twist map defined by a generating function:
 $$h(x,x')=\ds\frac{1}{2}(x-x')^2+\epsilon V(x), \quad  V(x)=\sum_{k=1}^\ell \ds\frac{a_k}{(2\pi k)^2}\cos(2\pi k x), \quad a_k\in\R,\quad \ell\geq 2,$$
 can exhibit in some range of parameters two $(p,q)$--minimizing orbits. Such maps are reversible
 i.e. there exists an involution $R$ such that  $f^{-1}=R\circ f\circ R$. $f$ factorizes as $f=I\circ R$, where $I$ is also an involution. The fixed point sets
 of the two involutions consist in two components. $\mbox{Fix}(R)=\Gamma_0\cup \Gamma'_0$, $\mbox{Fix}(I)=\Gamma_1\cup \Gamma'_1$, where:
 $$
 \Gamma_0: x=0,\quad \Gamma'_0: x=1/2, \quad \Gamma_1: x=y/2\:(\mbox{mod}\: 1), \quad \Gamma'_1: x=(y-1)/2\: (\mbox{mod}\: 1)$$

 Symmetric  $(p,q)$--periodic orbits (i.e. orbits that intersect $\mbox{Fix}(R)$) of a reversible twist map can undergo a Rimmer bifurcation \cite{mackay93}.
 Namely, if for fixed $a_k$, $k=\overline{1,\ell}$,  at some threshold $\epsilon=\epsilon_r$, one of the two  symmetric $(p,q)$--orbit changes its   extremal type either from  minimizing to minimaximing or conversely, and  two asymmetric $(p,q)$--orbits are born, then for $\epsilon$ in an interval $(\epsilon_r, \epsilon_c)$ the corresponding twist map has  two minimizing  and two minimaximizing $(p,q)$--orbits.
  To each of the two $(p,q)$--minimizing orbits starting at $(x_0, y_0)$, respectively at $(x'_0, y'_0)$, corresponds a $(p,q)$--minimizer,  $\x=(x_0, x_1, \ldots, x_{q-1})$, respectively  $\x'=(x'_0,x'_1, \ldots, x'_{q-1})\in\mathcal{X}_{pq}$.
The set of $(p,q)$--minimizers is in this case the union of two completely ordered subsets of $\mathcal{X}_{pq}$,
   $M_{pq}=\{\tau_{ij}(\x), i, j\in \Z\}\cup \{\tau_{k\ell}(\x'), k,\ell\in\Z\}$.  The elements in the second subset interlace those of the first subset.  Thus in $\mathcal{X}_{pq}/\tau_{01}$ there exist
  $2q$ minimizers. Between two consecutive $(p,q)$--minimizers there exists a $(p,q)$--minimaximizer.\par
\begin{figure}[h]
\begin{center}
\includegraphics{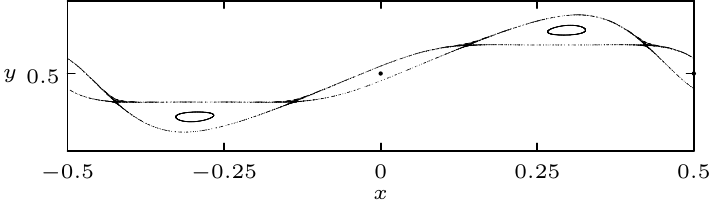}
\caption{\label{dynstdmap3hRimm}  The four $(1,2)$--periodic orbits of a three-harmonic standard map after a Rimmer bifurcation.}
\end{center}
\end{figure}
\begin{figure}[h]
\begin{center}
\includegraphics{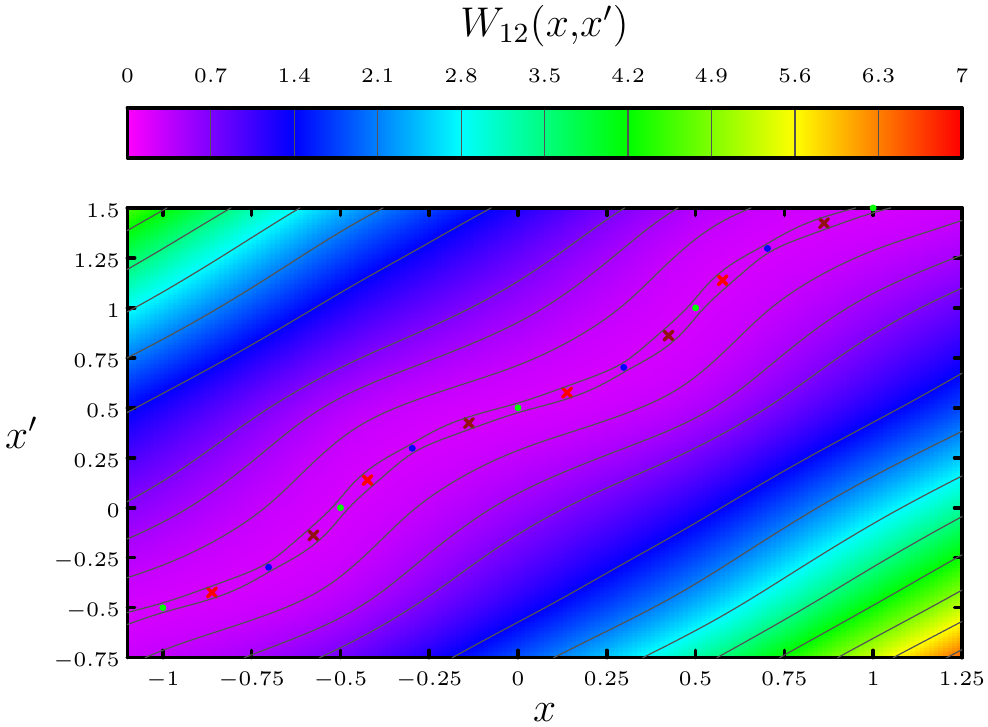}
\caption{\label{stdsin3W12Rimm}  The contour plot of the action $W_{12}$ associated to a generating function defined by a three harmonic potential. The associated twist map has four $(1,2)$--minimizers and four $(1,2)$--minimaximizers within $\mathcal{X}_{12}/\tau_{01}$.}
\end{center}
\end{figure}

  Such a case is shown in   Fig.\ref{stdsin3W12Rimm}, which  illustrates the contour plot of the action $W_{12}$ defined by the generating function corresponding to $\epsilon=1.2$, and the three-harmonic potential:
$$V(x)= \sum_{k=1}^3 \ds\frac{a_k}{(2k\pi)^2}\cos(2k\pi x),\quad   a_1=1, a_2=-0.3, a_3=0.2$$
The corresponding twist map exhibits two $(1,2)$--minimizing asymmetric orbits, respectively two $(1,2)$--minimaximizing symmetric orbits (Fig.\ref{dynstdmap3hRimm}).

The red, respectively  dark red crosses in Fig.\ref{stdsin3W12Rimm}, represent distinct $(1,2)$--minimizers $\x, \x'$, and a few of their integer translates, while blue and lightblue points are the interlacing distinct $(1,2)$--minimaximizers.

The asymmetric $(p,q)$--configurations present interest  and are also studied in condensed--matter physics  \cite{sasaki}.\par

From Proposition \ref{nostatincon} it follows that if $\x,\x'$ are two consecutive $(p,q)$--minimizers and $\y$ is a minimaximizer, $\x\prec\y\prec \x'$, then
any  stationary state, $\bf z$, ordered with respect to $M_{pq}$, and different from $\y$ can lie only within $[[\x, \x']]\setminus \left([\x,\y]\cup [\y, \x']\right)$.  \par
In many theoretical and numerical investigations of the standard map and the two--harmonic standard map reported in the literature, no ordered $(p,q)$--orbit different from the minimizing and minimaximizing orbit was revealed. This is due to the fact that the main tool used to decide the stability type of a $(p,q)$--periodic orbit is its residue. In \cite{petrisor} it was pointed out through a few examples that residue can be a misleading quantity in the characterization of a periodic orbit. It allows only to deduce the stability type of a periodic orbit (regular hyperbolic, elliptic or inverse hyperbolic) but not the right  extremal type of the $(p,q)$--stationary state corresponding to that orbit. Instead the twist number of a periodic orbit is much more relevant.  The twist number measures the average  rotation of tangent vectors under the action of the derivative of the twist map along a periodic orbit \cite{petrisor}.\par

In the following we show that the three harmonic standard map from the Example 4 \cite{petrisor} has a $(1,2)$--periodic orbit starting at at a point  $(0, 0.5)$, whose corresponding
$(1,2)$-configuration ${\bf z}=(z_0, z_1)$ is ordered with respect to $M_{12}$ and it is not a minimaximizer.
\par
\begin{figure}[h]
\begin{center}
\includegraphics{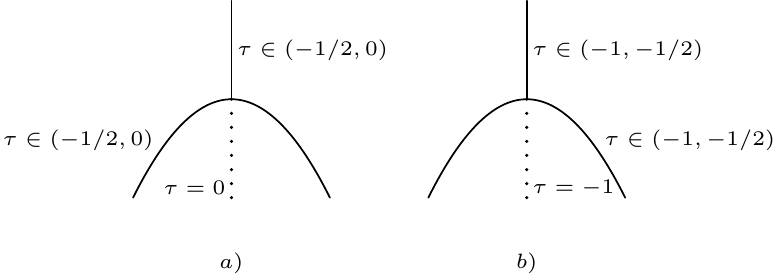}
\caption{\label{RimmerIllust}   Rimmer bifurcation that leads to the birth of two asymmetric orbits. In a) is the classical Rimmer bifurcation, and in b) a Rimmer type--bifurcation in which are involved orbits of large absolute twist.}
\end{center}
\end{figure}

The generating function that defines the map of interest  has the potential:
$$V(x)=-\ds\frac{0.18}{2\pi}\cos(2\pi x)+\ds\frac{0.42}{4\pi}\cos(4 \pi x)+\ds\frac{0.11}{6\pi}\cos(6\pi x)$$
The associated three harmonic standard map has four ordered $(1,2)$--orbits: a minimizing  one which is symmetric of twist number $\tau=0$, two asymmetric minimaximizing orbits
of twist number $\tau=-1/2$ (they are inverse hyperbolic orbits) and a symmetric orbit of twist number $-1$ (for details see \cite{petrisor}).
Connecting this three harmonic standard map with the integrable map one has deduced in \cite{petrisor}  that the minimaximizing orbit intersecting $\Gamma_0$
undergoes a sequence of bifurcations that leads to the decrease of its twist number. The last bifurcation is a  Rimmer--type  bifurcation. At the bifurcation threshold it turns from an orbit of $\tau\in(-1, -1/2)$ into an orbit of
$\tau=-1$, and two asymmetric orbits of $\tau=(-1,-1/2)$ are born.  This means that one of the two classical scenarios  of Rimmer bifurcation  illustrated  in Fig.\ref{RimmerIllust} a) occurs in our case as in  Fig.\ref{RimmerIllust} b).  Moreover the asymmetric orbits bifurcates further from $\tau\in(-1, -1/2)$ to $\tau=-1/2$. Each bifurcation is a bifurcation of an ordered orbit and thus after each threshold the old and the new born orbits are also ordered.\par
\begin{figure}[h]
\begin{center}
\includegraphics{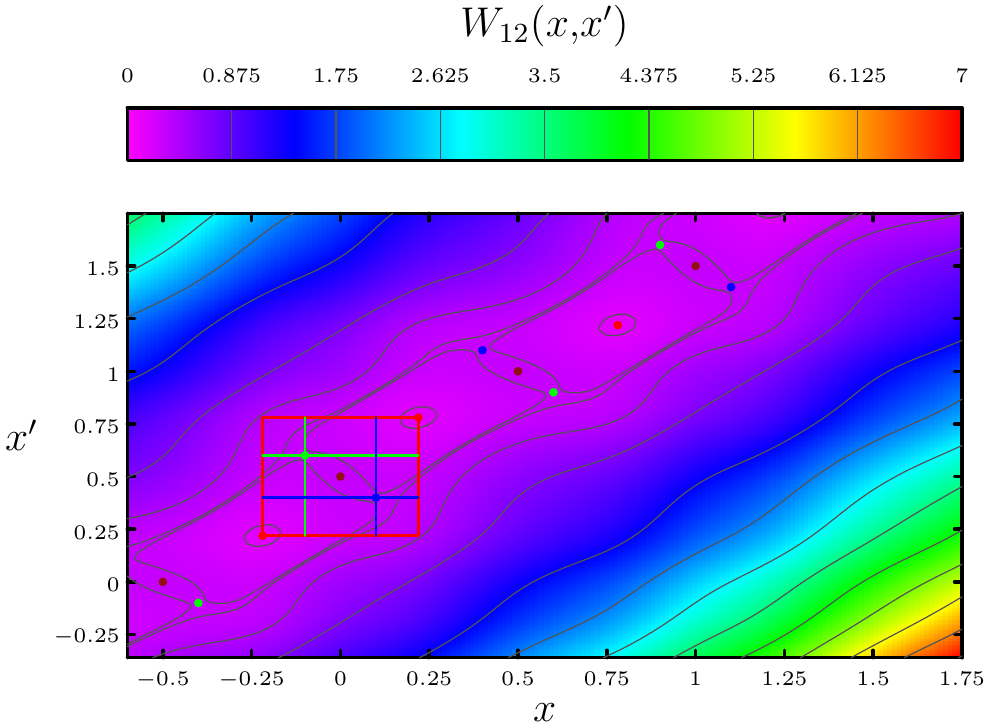}
\caption{\label{stdsinW12N}  The contour plot of the action $W_{12}$ associated to a generating function defined by a three harmonic potential. The corresponding twist map exhibits an ordered $(1,2)$--periodic orbit different from the minimizing and minimaximizing $(1,2)$--orbits.}
\end{center}
\end{figure}

The corresponding $(1,2)$--configurations are shown in Fig.\ref{stdsinW12N}.  To the  symmetric orbit of twist number $\tau=-1$ corresponds the dark red $(1,2)$--configurations. One can see that such a configuration that lies in an interval $[[\x, \x']]$ having as ends two consecutive  $(1,2)$--minimizers (the red points),  is incomparable with the $(1,2)$--minimaximizers  (colored in blue and green).

We note that the existence of such  unusual ordered $(p,q)$--configurations is ensured by the property  of the analyzed three-harmonic standard map, namely
it is a two--component strong folding region map \cite{petrisor}. \par\medskip
 All properties deduced in this note are illustrated in the space $\mathcal{X}_{12}$,  because it has  small dimension and the stationary states can be easily visualized. \par
The space of $(1,2)$--configurations is somewhat particular because if $\x=(x_0, x_1)$ is a stationary configuration, then both semilines of  the lines $X_0=x_0$,  $X_1=x_1$ (in the system of coordinates $(X_0, X_1)$), starting from $\x$ are boundaries for one of the order cones $V_{\pm}(\x)$. For $q>2$ this is not the case, and only a part from a hyperplan $X_i=x_i$, $i=\overline{0,q-1}$ is a boundary for an order cone.\par
More precisely  if $\x=(x_0, x_1, \ldots, x_{q-1})$ is a stationary state then the boundary of the cones $V_{\pm}(\x)$ is the union of the sets: $$B_i(\pm)=\{\z=(z_0, \ldots, z_{i-1}, x_i, z_{i+1}, \ldots, z_{q-1})\:|\:  z_j-x_j\geq 0,\:\, \mbox{for}\:\, V_+, z_j-x_j\leq 0 \:\, \mbox{for}\:\, V_-, j\neq i\}$$
$i=\ovl{0, q-1}$.

This remark is exploited in the following Proposition that  gives a sufficient condition ensuring that a $(p,q)$--periodic orbit  is unordered.
\bpr  Let  $(x_i, y_i)$, $i=\ovl{0,q-1}$, be   a $(p,q)$--minimizing orbit of a positive twist map $f$. If there is another $(p,q)$--periodic orbit that intersects one of the vertical lines $x=x_j$, $j=\ovl{0,q-1}$, then that orbit is unordered.\epr
Proof:  Let $\x'=(x'_0, x'_1, \ldots, x'_{q-1})$  be a $(p,q)$--minimizer in $\mathcal{X}_{pq}$ with $x'_0=x_j$, and  ${\bf z}=(z_0, z_1, \ldots, z_{q-1})$, $z_0=x_j$ a stationary state associated to the second orbit. $\bf z$ belongs to the hyperplane of equation $X_j=x_j$, but by Proposition \ref{nostateBcone}  ${\bf z}\notin V_{\pm}(\x')$.  Hence  $\bf z$ is  unordered.\par
\begin{cor}  If $f$ is a multi--harmonic standard map having a $(p,q)$--orbit that intersects the symmetry line $\Gamma_0$ or $\Gamma_1$ that is also intersected by
the $(p,q)$--minimizing orbit, then it is unordered.\end{cor}

\section{Conclusions} Using the duality between   an area preserving  twist map of the infinite cylinder, and the FK   model we deduced the location of stationary $(p,q)$--configurations in the space of all such configurations. We exploited the strict monotone property of the semiflow of the  action gradient in order to locate the regions where no equilibrium point of the gradient can lie.
Moreover we pointed out through an example that a twist map can exhibit ordered periodic orbits whose corresponding stationary configurations are neither action minimizing nor minimaximizing.  As far as we know no such an orbit was
identified before in the dynamics of a twist map.\par
The knowledge of forbidden subsets for stationary states as well as the subsets where ordered or unordered configurations can lie is useful
in numerical search for zeros of the action gradient, because the success of a  numerical method for detecting zeros of a vector field depends on the appropriate choice of the initial condition.  The existence of ordered stationary $(p,q)$--configurations,  whose associated periodic orbits have large absolute twist number can also be of  physical interest from the
point of view of the theory of commensurate--incommensurate phase transitions\par
The   numerical assessment of  the Boyland--Hall criterion \cite{leagemcK} led to the conclusion that it is not as successful as other converse KAM criteria, because  unordered periodic orbits of small period were detected only for large values of the perturbation parameter, for which other methods ensured the breakup of all invariant circles of the standard map. In our theoretical and numerical study of uni--component and two--component strong folding region twist maps,  started in \cite{petrisor}, it appears that just after the breakdown of an invariant circle, the nearby unordered orbits  have large period, not small. Thus one expects that this criterion can work  well if we look for unordered periodic orbits of large period.
   In order to confirm and illustrate this behaviour we need to know where in the space of $(p,q)$--configurations we must look for  unordered stationary configurations.
 \par

\par
 \end{document}